# TIME-FREQUENCY ANALYSIS AND HARMONIC GAUSSIAN FUNCTIONS

Tokiniaina Ranaivoson, Raoelina Andriambololona, Rakotoson Hanitriarivo

Theoretical Physics Department
Institut National des Sciences et Techniques Nucléaires (INSTN-Madagascar), Antananarivo, Madagascar

*instn@moov.mg, tokhiniaina@gmail.com (T. Ranaivoson), raoelinasp@yahoo.fr (Raoelina Andriambololona),*
*jacquelineraoelina@hotmail.com (Raoelina Andriambololona), infotsara@gmail.com (R. Hanitriarivo)*

***Abstract***: A method for time-frequency analysis is given. The approach utilizes properties of Gaussian distribution, properties of Hermite polynomials and Fourier analysis. We begin by the definitions of a set of functions called harmonic Gaussian functions. Then these functions are used to define a set of transformations, noted $\mathcal{T}_n$, which associate to a function $\psi$, of the time variable $t$, a set of functions $\Psi_n$ which depend on time, frequency and frequency (or time) standard deviation. Some properties of the transformations $\mathcal{T}_n$ and the functions $\Psi_n$ are given. It is proved in particular that the square of the modulus of each function $\Psi_n$ can be interpreted as a representation of the energy distribution of the signal, represented by the function $\psi$, in the time-frequency plane for a given value of the frequency (or time) standard deviation. It is also shown that the function $\psi$ can be recovered from the functions $\Psi_n$.

**Keywords**: Time-frequency analysis; Signal; Energy distribution; Gaussian distribution; Hermite polynomials

## 1. Introduction

The time-frequency analysis is an important field of study and has many important applications. In this analysis, a main purpose is to have a good representation of a signal and/or the distribution of its energy both in time and frequency. Because of the uncertainty relation, it is not possible to have a rigorous representation of the time-frequency distribution of the energy of a signal at any scale of resolution in time and frequency. Let $\psi$ be a function, of the time variable $t$, which represent a signal and let $\tilde{\psi}$ be the Fourier Transform of $\psi$

$$\tilde{\psi}(\omega) = \frac{1}{\sqrt{2\pi}} \int_{-\infty}^{+\infty} \psi(t)\, e^{-i\omega t} dt \quad (1.1a)$$

$$\psi(t) = \frac{1}{\sqrt{2\pi}} \int_{-\infty}^{+\infty} \tilde{\psi}(\omega)\, e^{i\omega t} d\omega \quad (1.1b)$$

The variable $\omega = 2\pi f$ is the angular frequency (if $f$ is the frequency). We will use $\omega$ as the frequency variable. The energy $E$ of the signal is defined by the relation

$$E = \int_{-\infty}^{+\infty} |\psi(t)|^2\, dt = \int_{-\infty}^{+\infty} |\tilde{\psi}(\omega)|^2\, d\omega \quad (1.2)$$

The mean values $\mu_t$ and $\mu_\omega$ of time and frequency associated to the function $\psi$ are

$$\mu_t = \frac{1}{E}\int t|\psi(t)|^2 dt \qquad \mu_\omega = \frac{1}{E}\int \omega|\tilde{\psi}(\omega)|^2 d\omega \quad (1.3)$$

The time and frequency standard deviations are

$$\sigma_t = \sqrt{\frac{1}{E}\int (t-\mu_t)^2 |\psi(t)|^2\, dt} \quad (1.4)$$

$$\sigma_\omega = \sqrt{\frac{1}{E}\int (\omega-\mu_\omega)^2 |\tilde{\psi}(\omega)|^2\, d\omega} \quad (1.5)$$

Using properties of Fourier transform and Cauchy-Schwarz inequality, it may be shown that we have the inequality

$$\sigma_t \sigma_\omega \geq \frac{1}{2} \quad (1.6)$$

The inequality (1.6) is the uncertainty relation between time and frequency.

Many attempts have been done to formulate methods to give a good solution to the problem of time frequency analysis. Most of those methods can be considered as based on the use of bilinear time-frequency distribution or short-time Fourier transform [1], [2], [3], [4], [5],[6]. Wavelets transforms can also considered as related to the time frequency analysis [7].

One of the well known bilinear distribution which can be considered is the Wigner-Ville distribution [8],[9], [10],[11]. In signal analysis, the Wigner-Ville distribution may be defined as the function of the two variables time $t$ and frequency $\omega$, associated to a function $\psi$ of the variable $t$ by the relation

$$W_\psi(t,\omega) = \frac{1}{2\pi}\int_{-\infty}^{+\infty} \psi^*\left(t-\frac{u}{2}\right)\psi\left(t+\frac{u}{2}\right) e^{-i\omega u} du$$

$$= \frac{1}{2\pi}\int_{-\infty}^{+\infty} \tilde{\psi}^*\left(\omega-\frac{\xi}{2}\right)\tilde{\psi}\left(\omega+\frac{\xi}{2}\right) e^{i\xi t} d\xi \quad (1.7)$$

where $\tilde{\psi}$ is the Fourier transform of $\psi$.
If $E$ is the energy of the signal represented by the function $\psi$ as defined in the relation (1.2), we have the relations

$$E = \int_{-\infty}^{+\infty} |\psi(t)|^2\, dt = \int_{-\infty}^{+\infty} |\tilde{\psi}(\omega)|^2\, d\omega \quad (1.8a)$$

$$= \int_{-\infty}^{+\infty} \int_{-\infty}^{+\infty} W_\psi(t,\omega)\, dt d\omega \qquad (1.8b)$$

The relations (1.8$a$) and (1.8$b$) suggest that the Wigner-Ville distribution may be interpreted as the density of the energy of a signal in the time-frequency plane. But it can be shown that the function $W_\psi$ as defined in the relation (1.7) is not positive definite. In fact, using Cauchy-Schwarz inequality in the space $\mathcal{L}^2$ of Lebesgue square-integrable functions and the relations (1.7) and (1.8a) we can deduce the inequality

$$|W_\psi(t,\omega)| \leq \frac{E}{\pi} \Leftrightarrow -\frac{E}{\pi} \leq W_\psi(t,\omega) \leq \frac{E}{\pi} \qquad (1.9)$$

This propriety of the Wigner-Ville distribution doesn't allow to interpret it easily as a good representation of the energy density. Depth studies on the Wigner-Ville distribution and its properties have been done by many authors to give solutions to this problem, as example of solution is the introduction of smoothing methods which reduce the amount of negativity.

As mentioned previously, another method used in the time frequency analysis is the short-time Fourier transform. For a signal represented by a function $\psi$ of the time variable $t$, a short-time Fourier transform $\Psi$ of $\psi$ may be defined by the relation

$$\Psi(T,\Omega) = \frac{1}{\sqrt{2\pi}} \int_{-\infty}^{+\infty} \psi(t)\, g^*(t-T) e^{-i\Omega t} dt \quad (1.10)$$

The function $g$ is the window function. $g(t)$ has significant values in the vicinity of $t=0$ and tends to zero outside this range.

We have the relation

$$\int_{-\infty}^{+\infty}\int_{-\infty}^{+\infty} |\Psi(T,\Omega)|^2\, dTd\Omega = \int_{-\infty}^{+\infty}\int_{-\infty}^{+\infty} |\psi(t) g(t-T)|^2 dT\, dt$$

so if the window function $g$ satisfies the relation

$$\int_{-\infty}^{+\infty} |g(t-T)|^2 dT = 1$$

we have

$$\int_{-\infty}^{+\infty}\int_{-\infty}^{+\infty} |\Psi(T,\Omega)|^2\, dTd\Omega = \int_{-\infty}^{+\infty} |\psi(t)|^2 dt = E \quad (1.11)$$

According to this relation, the square of the modulus $|\Psi(T,\Omega)|^2$ of the short time Fourier transform may be interpreted as a representation of the energy density of the signal in time-frequency plane. And unlike the case of the Wigner-Ville distribution, it is a positive definite distribution. A special case of short-time Fourier transform is the Gabor transform in which the window function $g$ is a gaussian function

$$g(t) = 2^{1/4} e^{-\pi(t)^2} \qquad (1.12)$$

If we introduce the Gabor function

$$\varphi(t,T,\Omega) = g(t-T) e^{i\Omega t} = 2^{1/4} e^{-\pi(t-T)^2} e^{i\Omega t} \quad (1.13)$$

a Gabor transform of a function $\psi$ may be defined by the relation

$$\Psi(T,\Omega) = \frac{1}{\sqrt{2\pi}} \int_{-\infty}^{+\infty} \varphi^*(t,T,\Omega)\psi(t)\, dt \qquad (1.14)$$

In this paper, we tackle the problem of time-frequency analysis with the introduction of a set of transformations, noted $\mathcal{T}_n$ ($n \in \mathbb{N}$), which associate to a function $\psi$ of the time variable $t$ a set of functions $\Psi_n$ which depend on time, frequency and frequency (or time) standard deviation. We show that the square of the modulus of each function $\Psi_n$ may be interpreted as a representation of the time-frequency distribution of the energy of the signal corresponding to the function $\psi$ (relations (4.9) and (4.10)). Then we show that the function $\psi$ can be recovered from the functions $\Psi_n$ (section 5).

For the definitions of the transformations $\mathcal{T}_n$, we introduce a set of functions $\{\varphi_n\}$ that we will call harmonic Gaussian functions. In a certain point of view, the functions $\varphi_n$ may be seen as generalizations of Gabor functions and the transformations $\mathcal{T}_n$ as generalizations of Gabor transformations. However, our results show obviously that there are differences between our method and the Gabor's one. Compared to the Wigner-Ville distribution, the energy distributions that we introduce in our method has the advantage to be positive definite.

## 2. Harmonic Gaussian Functions

For positive integers $n$, let us define a set of orthonormalized functions denoted $\varphi_n$ such as

$$\varphi_n(t,T,\Omega,\Delta t) = \frac{H_n\left(\frac{t-T}{\sqrt{2}\Delta t}\right)}{\sqrt{2^n n! \sqrt{2\pi}\Delta t}} e^{-\left(\frac{t-T}{2\Delta t}\right)^2 + i\Omega t} \qquad (2.1)$$

$$\int_{-\infty}^{+\infty} \varphi_m^*(t,T,\Omega,\Delta t)\varphi_n(t,T,\Omega,\Delta t) dt = \delta_{nm} \qquad (2.2)$$

in which $H_n$ is Hermite polynomial of degree $n$. Useful properties of Hermite polynomials are recalled in the appendix A. The set of functions $\{\varphi_n\}_{n\in\mathbb{N}}$ is an orthonormal basis of the vectorial space $\mathcal{L}^2$ of Lebesgue square integrable functions.

By taking into account the relations (A.6) and (A.7) in the appendix A, we may establish easily the relations

$$\int_{-\infty}^{+\infty} t|\varphi_n(t,T,\Omega,\Delta t)|^2\, dt = T \qquad (2.3)$$

$$\int_{-\infty}^{+\infty} (t-T)^2 |\varphi_n(t,T,\Omega,\Delta t)|^2\, dt = (2n+1)(\Delta t)^2 \quad (2.4)$$

According to the relations (2.3) and (2.4), we call $T$ the temporal mean value, $(\sigma_t)^2 = (\Delta t_n)^2 = (2n+1)(\Delta t)^2$ the time variance and $\sigma_t = \Delta t_n$ the time standard deviation corresponding to $\varphi_n$.

According to the result given in the appendix B, the expression of the Fourier transform $\tilde{\varphi}_n$ of the function $\varphi_n$ is

$$\tilde{\varphi}_n(\omega, T, \Omega, \Delta\omega) = \frac{i^n H_n\left(\frac{\omega - \Omega}{\sqrt{2}\Delta\omega}\right)}{\sqrt{2^n n! \sqrt{2\pi}\Delta\omega}} e^{-\left(\frac{\omega - \Omega}{2\Delta\omega}\right)^2 - iT(\omega - \Omega)} \quad (2.5)$$

in which $\Delta\omega$ is related to $\Delta t$ by the relation

$$\Delta t \Delta\omega = \frac{1}{2} \quad (2.6)$$

From now on, we assume that the relation (2.6) is always fulfilled by $\Delta t$ and $\Delta\omega$.

As for the case of $\varphi_n$, we may establish in the case of $\tilde{\varphi}_n$ the relations

$$\int_{-\infty}^{+\infty} \tilde{\varphi}_m^*(\omega, T, \Omega, \Delta\omega) \tilde{\varphi}_n(\omega, T, \Omega, \Delta\omega) d\omega = \delta_{nm} \quad (2.7)$$

$$\int_{-\infty}^{+\infty} \omega |\tilde{\varphi}_n(\omega, T, \Omega, \Delta\omega)|^2 d\omega = \Omega \quad (2.8)$$

$$\int_{-\infty}^{+\infty} (\omega - \Omega)^2 |\tilde{\varphi}_n(\omega, T, \Omega, \Delta\omega)|^2 d\omega = (2n+1)(\Delta\omega)^2 \quad (2.9)$$

According to the relations (2.8) and (2.9), we call $\Omega$ the frequency mean value, $(\sigma_\omega)^2 = (\Delta\omega_n)^2 = (2n+1)(\Delta\omega)^2$ the frequency variance and $\sigma_\omega = \Delta\omega_n$ the frequency standard deviation corresponding to $\tilde{\varphi}_n$.

Because of the similarity between the expression (2.1) of a function $\varphi_n$ and the expression of wave functions of a linear harmonic oscillator in quantum mechanics [12], we will call a function $\varphi_n$ a "Harmonic Gaussian Function of degree $n$". According to the above results, a harmonic Gaussian function is characterized by its time mean value $T$, its frequency mean value $\Omega$, its time variance $(\Delta t_n)^2 = (2n+1)(\Delta t)^2$ and its frequency variance $(\Delta\omega_n)^2 = (2n+1)(\Delta\omega)^2$. $\Delta t$ and $\Delta\omega$ are related by the relation (2.6).

## 3. Representations of a signal with functions defined in the time-frequency plane

Let $\psi$ be a function of the time variable $t$ which represent a signal, $\psi$ is an element of the space $\mathcal{L}^2$ of Lebesgue square-integrable functions. Let be $\tilde{\psi}$ the Fourier transform of $\psi$, $\tilde{\psi}$ is an element of $\mathcal{L}^2$ too. We denote $\mathbb{F}$ the functions space generated by functions $\Psi$ defined in the time-frequency plane i.e. functions of the two riables $T, \Omega$ and which is Lebesgue square-integrable in the time-frequency plane $\{(T, \Omega)\}$

$$\int_{-\infty}^{+\infty} \int_{-\infty}^{+\infty} |\Psi(T, \Omega)|^2 dT d\Omega = K < +\infty \quad (3.1)$$

For a positive integer $n$, let be $\mathcal{T}_n$ the application of $\mathcal{L}^2$ to $\mathbb{F}$ defined by

$$\mathcal{T}_n: \quad \mathcal{L}^2 \to \mathbb{F}$$
$$\psi \mapsto \Psi_n = \mathcal{T}_n(\psi)$$

$$\Psi_n(T, \Omega, \Delta\omega) = \frac{1}{\sqrt{2\pi}} \int_{-\infty}^{+\infty} \varphi_n^*(t, T, \Omega, \Delta t) \psi(t) dt \quad (3.2)$$

$$= \frac{1}{\sqrt{2\pi}} \int_{-\infty}^{+\infty} \tilde{\varphi}_n^*(\omega, T, \Omega, \Delta\omega) \tilde{\psi}(\omega) d\omega \quad (3.3)$$

In the relations (3.2) and (3.3), the function $\varphi_n^*$ is the complex conjugate of the harmonic Gaussian function $\varphi_n$ defined in the relation (2.1) and the function $\tilde{\varphi}_n^*$ is the complex conjugate of the Fourier transform $\tilde{\varphi}_n$ (relation 2.5) of the function $\varphi_n$.

**Theorem 1**

The transformation $\mathcal{T}_n$ is linear

$$\forall \psi, \phi \in \mathcal{L}^2, \quad \forall \lambda, \mu \in \mathbb{C}:$$
$$\mathcal{T}_n(\lambda\psi + \mu\phi) = \lambda\mathcal{T}_n(\psi) + \mu\mathcal{T}_n(\phi) \quad (3.4)$$

**Proof**

Using the relation (3.2), we have

$$\mathcal{T}_n[(\lambda\psi + \mu\phi)(t)]$$
$$= \frac{1}{\sqrt{2\pi}} \int_{-\infty}^{+\infty} \varphi_n^*(t, T, \Omega, \Delta t)[\lambda\psi(t) + \mu\phi(t)] dt$$
$$= \frac{\lambda}{\sqrt{2\pi}} \int_{-\infty}^{+\infty} \varphi_n^*(t, T, \Omega, \Delta t) \psi(t) dt$$
$$+ \frac{\mu}{\sqrt{2\pi}} \int_{-\infty}^{+\infty} \varphi_n^*(t, T, \Omega, \Delta t) \phi(t) dt$$
$$= \lambda\mathcal{T}_n[\psi(t)] + \mu\mathcal{T}_n[\phi(t)]$$
$$\Rightarrow \mathcal{T}_n(\lambda\psi + \mu\phi) = \lambda\mathcal{T}_n(\psi) + \mu\mathcal{T}_n(\phi)$$

According to the definitions of the transformations $\mathcal{T}_n$, a function $\Psi_n$ may be considered as a time-frequency representation of the signal corresponding to the function $\psi$ for a given time resolution and frequency resolution characterized by the time standard deviation (or by the frequency standard deviation) corresponding to the harmonic Gaussian functions which is used for the definition of the transformation $\mathcal{T}_n$. The next section gives more justification to this point of view.

## 4. Representations of the energy distribution in the time-frequency plane

We consider a signal represented by the functions $\psi, \tilde{\psi}$ and $\Psi_n$ as defined in the previous section. The energy $E$ of the signal is given by the relations

$$E = \int_{-\infty}^{+\infty} |\psi(t)|^2 dt = \int_{-\infty}^{+\infty} |\tilde{\psi}(\omega)|^2 d\omega \quad (4.1)$$

In this section, we prove that the square of the modulus of the function $\Psi_n$ for a given $n \in \mathbb{N}$ may be interpreted as a representation of the distribution of the energy of the signal in the time frequency plane at a given scale of resolution characterized by frequency (or time) standard deviation.

*Theorem 2*

Let us consider the functions

$$p_n(T, \Delta t) = \int_{-\infty}^{+\infty} |\Psi_n(T, \Omega, \Delta\omega)|^2 \, d\Omega \quad (4.2)$$

$$\rho_n(\Omega, \Delta\omega) = \int_{-\infty}^{+\infty} |\Psi_n(T, \Omega, \Delta\omega)|^2 \, dT \quad (4.3)$$

We have the relations

$$p_n(T, \Delta t) = \int_{-\infty}^{+\infty} |\varphi_n(t, T, \Omega, \Delta t)|^2 |\psi(t)|^2 \, dt \quad (4.4a)$$

$$= \int_{-\infty}^{+\infty} \frac{\left|H_n\left(\frac{t-T}{\sqrt{2}\Delta t}\right)\right|^2}{2^n n! \sqrt{2\pi}\Delta t} e^{-\left(\frac{t-T}{\sqrt{2}\Delta t}\right)^2} |\psi(t)|^2 \, dt \quad (4.4b)$$

$$\rho_n(\Omega, \Delta\omega) = \int_{-\infty}^{+\infty} |\tilde{\varphi}_n(\omega, T, \Omega, \Delta\omega)|^2 |\tilde{\psi}(\omega)|^2 \, d\omega \quad (4.5a)$$

$$= \int_{-\infty}^{+\infty} \frac{\left|H_n\left(\frac{\omega-\Omega}{\sqrt{2}\Delta\omega}\right)\right|^2}{2^n n! \sqrt{2\pi}\Delta\omega} e^{-\left(\frac{\omega-\Omega}{\sqrt{2}\Delta\omega}\right)^2} |\tilde\psi(\omega)|^2 d\omega \quad (4.5b)$$

**Proof**

To prove the relations $(4.4a)$ and $(4.4b)$, we use the relation $(3.2)$, then we have

$$p_n(T, \Delta t) = \int_{-\infty}^{+\infty} |\Psi_n(T, \Omega, \Delta t)|^2 \, d\Omega$$

$$= \frac{1}{2\pi}[\int_{-\infty}^{+\infty}\int_{-\infty}^{+\infty}\int_{-\infty}^{+\infty} \varphi_n(t,T,\Omega,\Delta t)\varphi_n^*(t',T,\Omega,\Delta t)$$

$$\psi(t)\psi^*(t')]\,dt\,dt'\,d\Omega$$

$$= \{\int_{-\infty}^{+\infty}\int_{-\infty}^{+\infty}[\frac{1}{2\pi}\int_{-\infty}^{+\infty}\varphi_n(t,T,\Omega,\Delta t)\varphi_n^*(t',T,\Omega,\Delta t)d\Omega]$$

$$\psi(t)\psi^*(t')\}\,dt\,dt'$$

We have for the $\Omega$ integration

$$\frac{1}{2\pi}\int_{-\infty}^{+\infty}\varphi_n(t,T,\Omega,\Delta t)\varphi_n^*(t',T,\Omega,\Delta t)d\Omega$$

$$= \int_{-\infty}^{+\infty} \frac{H_n\left(\frac{t-T}{\sqrt{2}\Delta t}\right)H_n\left(\frac{t'-T}{\sqrt{2}\Delta t}\right)}{2^n n! \sqrt{2\pi}\Delta t} e^{-\left(\frac{t-T}{2\Delta t}\right)^2+\left(\frac{t'-T}{2\Delta t}\right)^2} \frac{e^{i\Omega(t-t')}}{2\pi} d\Omega$$

$$= \frac{H_n\left(\frac{t-T}{\sqrt{2}\Delta t}\right)H_n\left(\frac{t'-T}{\sqrt{2}\Delta t}\right)}{2^n n! \sqrt{2\pi}\Delta t} e^{-\left(\frac{t-T}{2\Delta t}\right)^2+\left(\frac{t'-T}{2\Delta t}\right)^2} \delta(t-t')$$

in which $\delta(t-t')$ is the Dirac's distribution

$$\Rightarrow p_n(T, \Delta t)$$

$$= \int_{-\infty}^{+\infty}\int_{-\infty}^{+\infty}[\frac{H_n\left(\frac{t-T}{\sqrt{2}\Delta t}\right)H_n\left(\frac{t'-T}{\sqrt{2}\Delta t}\right)}{2^n n! \sqrt{2\pi}\Delta t} e^{-\left(\frac{t-T}{2\Delta t}\right)^2+\left(\frac{t'-T}{2\Delta t}\right)^2}$$

$$\psi(t)\psi^*(t')\delta(t-t')]dtdt'$$

$$= \int_{-\infty}^{+\infty} \frac{\left|H_n\left(\frac{t-T}{\sqrt{2}\Delta t}\right)\right|^2}{2^n n! \sqrt{2\pi}\Delta t} e^{-\frac{1}{2}\left(\frac{t-T}{\Delta t}\right)^2} |\psi(t)|^2 dt$$

$$= \int_{-\infty}^{+\infty} |\varphi_n(t, T, \Omega, \Delta t)|^2 |\psi(t)|^2 \, dt$$

The proof of the relations $(4.5a)$ and $(4.5b)$ may be obtained by analogy.

The relations $(4.2)$ and $(4.3)$ lead to the relations

$$\int_{-\infty}^{+\infty}\int_{-\infty}^{+\infty} |\Psi_n(T, \Omega, \Delta\omega)|^2 \, dTd\Omega = \int_{-\infty}^{+\infty} p_n(T, \Delta t) \, dT$$

$$= \int_{-\infty}^{+\infty} \rho_n(\Omega, \Delta\omega) \, d\Omega \quad (4.6)$$

From the relations $(4.4b)$ and $(4.5b)$ we can also deduced the relations

$$\lim_{\Delta t \to 0} p_n(T, \Delta t) = |\psi(T)|^2 \quad (4.7)$$

$$\lim_{\Delta\omega \to 0} \rho_n(\Omega, \Delta\omega) = |\tilde\psi(\Omega)|^2 \quad (4.8)$$

But the two limits are not to be performed simultaneously because of the relation $(2.6)$

**Theorem 3**

We have the relations

$$\int_{-\infty}^{+\infty}\int_{-\infty}^{+\infty} |\Psi_n(T, \Omega, \Delta\omega)|^2 \, dTd\Omega$$

$$= \int_{-\infty}^{+\infty} p_n(T, \Delta t) \, dT = \int_{-\infty}^{+\infty} |\psi(t)|^2 \, dt = E \quad (4.9)$$

$$= \int_{-\infty}^{+\infty} \rho_n(\Omega, \Delta\omega) \, d\Omega = \int_{-\infty}^{+\infty} |\tilde\psi(\omega)|^2 \, d\omega = E \quad (4.10)$$

**Proof**

To prove the relation $(4.9)$, we utilize the relations $(4.6)$ and $(4.4b)$

$$\int_{-\infty}^{+\infty}\int_{-\infty}^{+\infty} |\Psi_n(T, \Omega, \Delta\omega)|^2 \, dTd\Omega = \int_{-\infty}^{+\infty} p_n(T, \Delta t) \, dT$$

$$= \int_{-\infty}^{+\infty}\int_{-\infty}^{+\infty} \frac{\left|H_n\left(\frac{t-T}{\sqrt{2}\Delta t}\right)\right|^2}{2^n n! \sqrt{2\pi}\Delta t} e^{-\left(\frac{t-T}{\sqrt{2}\Delta t}\right)^2} |\psi(t)|^2 dt dT$$

$$= \int_{-\infty}^{+\infty} \left[ \int_{-\infty}^{+\infty} \frac{\left| H_n\left(\frac{t-T}{\sqrt{2}\Delta t}\right) \right|^2}{2^n n! \sqrt{2\pi}\Delta t} e^{-\left(\frac{t-T}{\sqrt{2}\Delta t}\right)^2} dT \right] |\psi(t)|^2 \, dt$$

$$= \int_{-\infty}^{+\infty} |\psi(t)|^2 \, dt = E$$

The relations (4.10) may be obtained by analogy from the relations (4.6) and (4.5b) instead of the relation (4.4b).

The analysis of the results in the theorems 2 and 3 leads to the following remarks and interpretations:
- The function $p_n(T, \Delta t)$ may be interpreted as a representation of the instantaneous power distribution of the signal at the scale of resolution characterized by the time standard deviation.
- The function $\rho_n(\Omega, \Delta \omega)$ may be interpreted as a representation of the spectral energy distribution of the signal at the scale of resolution characterized by the frequency standard deviation.
- The function $|\Psi_n(T, \Omega, \Delta \omega)|^2$ may be interpreted as a representation of the energy distribution of the signal in the time-frequency plane at the scale of resolution characterized by the time (or frequency) standard deviation.

Because of the existence of the parameter $n$, we may call $p_n(T, \Delta t)$ a representation of the instantaneous power distribution of the signal at order $n$, $\rho_n(\Omega, \Delta \omega)$ a representation of the spectral energy distribution of the signal at order $n$ and $|\Psi_n(T, \Omega, \Delta \omega)|^2$ a representation of the time-frequency distribution of the energy of the signal at order $n$.

## 5. Recovering of the original function and another expression for the energy

We can recover the function $\psi$ from the functions $\Psi_n$.

*Theorem 4*

We have the relations

$$\psi(t) = \sum_n \sqrt{2\pi} \Psi_n(T, \Omega, \Delta \omega) \, \varphi_n(t, T, \Omega, \Delta t) \quad (5.1)$$

$$E = \int_{-\infty}^{+\infty} |\psi(t)|^2 \, dt = \sum_n 2\pi |\Psi_n(T, \Omega, \Delta \omega)|^2 \quad (5.2)$$

**Proof**

Let us expand the function $\psi$ in the basis $\{\varphi_n\}_{n \in \mathbb{N}}$

$$\psi(t) = \sum_n C_n(T, \Omega, \Delta t) \, \varphi_n(t, T, \Omega, \Delta t)$$

$C_n(T, \Omega, \Delta t)$ are the components of $\psi$ in the basis $\{\varphi_n\}_{n \in \mathbb{N}}$. According to the relation (2.2), the basis $\{\varphi_n\}_{n \in \mathbb{N}}$ is an orthonormal basis, so we have

$$C_n(T, \Omega, \Delta t) = \int_{-\infty}^{+\infty} \varphi_n^*(t, T, \Omega, \Delta t) \psi(t) \, dt$$

$$= \sqrt{2\pi} \, \Psi_n(T, \Omega, \Delta \omega)$$

The last equality follows from the relation (3.2).
The relation (5.2) can be deduced easily from the relation (5.1) and the orthonormality of the basis $\{\varphi_n\}_{n \in \mathbb{N}}$.

*Theorem 5*

For the recovering of the function $\psi$, we have also the relation

$$\psi(t) = \frac{1}{\sqrt{2\pi}} \int_{-\infty}^{+\infty} \int_{-\infty}^{+\infty} \Psi_n(T, \Omega, \Delta \omega) \, \varphi_n(t, T, \Omega, \Delta t) \, dT d\Omega$$

**Proof**

$$\frac{1}{\sqrt{2\pi}} \int_{-\infty}^{+\infty} \int_{-\infty}^{+\infty} \Psi_n(T, \Omega, \Delta \omega) \, \varphi_n(t, T, \Omega, \Delta t) \, dT d\Omega =$$

$$\frac{1}{2\pi} \int_{-\infty}^{+\infty} \int_{-\infty}^{+\infty} \int_{-\infty}^{+\infty} \varphi_n^*(t', T, \Omega, \Delta t) \, \psi(t') \varphi_n(t, T, \Omega, \Delta t) dt' dT d\Omega$$

$$\int_{-\infty}^{+\infty} \int_{-\infty}^{+\infty} \int_{-\infty}^{+\infty} \frac{H_n\left(\frac{t-T}{\sqrt{2}\Delta t}\right) H_n\left(\frac{t'-T}{\sqrt{2}\Delta t}\right)}{2^n n! \sqrt{2\pi}\Delta t} \psi(t') e^{-\left[\left(\frac{t-T}{2\Delta t}\right)^2 + \left(\frac{t'-T}{2\Delta t}\right)^2\right]}$$
$$\frac{e^{i\omega(t-t')}}{2\pi} d\Omega dt' dT$$

$$= \int_{-\infty}^{+\infty} \int_{-\infty}^{+\infty} \frac{H_n\left(\frac{t-T}{\sqrt{2}\Delta t}\right) H_n\left(\frac{t'-T}{\sqrt{2}\Delta t}\right)}{2^n n! \sqrt{2\pi}\Delta t} \psi(t') e^{-\left[\left(\frac{t-T}{2\Delta t}\right)^2 + \left(\frac{t'-T}{2\Delta t}\right)^2\right]}$$
$$\delta(t-t') dt' dT$$

$$= \int_{-\infty}^{+\infty} \frac{\left| H_n\left(\frac{t-T}{\sqrt{2}\Delta t}\right) \right|^2}{2^n n! \sqrt{2\pi}\Delta t} \psi(t) e^{-\left[\left(\frac{t-T}{\sqrt{2}\Delta t}\right)^2\right]} dT$$

$$= \psi(t)$$

## 6. Conclusions

We may conclude that the introduction of harmonic Gaussian functions gives a possibility to establish new methods to tackle the problem of time-frequency analysis. Our approach differs particularly with other ones in the fact that we introduce simultaneously a set of functions which allows to have a set of possible representations of the energy distribution according to the values of the parameter $n$ in the expression of the harmonic Gaussian function $\varphi_n$. And it is possible to have more possible representations according to the values of the frequency (or time) standard deviation.

As shown in the relations (4.9) and (4.10), our approach allows also to have, with the set of the possible representations of the energy distribution in time-frequency plane, the sets of possible representations of the instantaneous power distribution (functions $p_n$) and possible representations of the spectral energy distribution (functions $\rho_n$).

Another interesting result is also the possibility of recovering the original function $\psi$ which represents the signal from the functions $\Psi_n$ (section 5). This recovering may be obtained by making a summation on the index $n$ or by making an integration on the time-frequency plane $\mathcal{P} = \{(T, \Omega)\}$. The first expansion needs the knowledge of the values of all the functions $\Psi_n$ for all $n$ for given values of $T, \Omega$ and $\Delta\omega$. The second one needs the knowledge of the expression of one function $\Psi_n$ for given values of $n$ and $\Delta\omega$.

The expression of the energy of the signal given in the relation (5.2), which is a direct consequence of the relation (5.1), may allows also to interpret $2\pi|\Psi_n(T, \Omega, \Delta\omega)|^2$ as the part of the energy of the signal which corresponds to the harmonic Gaussian function of degree $n$ for given values of $T, \Omega$ and $\Delta\omega$.

More depth studies on the physical meaning of the obtained mathematical results may give more interesting and useful insights for the theory of signal analysis and signal processing.

# Appendix A

**Useful properties of Hermite Polynomials [12]**

For positive integers $n$, the Hermite polynomial $H_n(x)$ can be defined by the relation

$$H_n(x) = (-1)^n e^{x^2} \frac{d^n}{dx^n}\left(e^{-x^2}\right) \quad (A.1)$$

The following properties may be established.

**Property 1**

For any $n \in \mathbb{N}$ and for any $x \in \mathbb{R}$, we have the recurrence relation:

$$H_{n+1}(x) = 2xH_n(x) - H_n'(x) \quad (A.2)$$

**Property 2**

Let us expand the Hermite polynomials

$$H_n(x) = \sum_{k=0}^{n} a_k^n x^k$$

Then we have the following relations:
$a_n^n = 2^n$ (coefficient of $x^n$ in $H_n(x)$)
$a_{n-1}^n = 0$ (coefficient of $x^{n-1}$ in $H_n(x)$)
$a_k^n = 2a_{k-1}^{n-1} - (p+1)a_{k+1}^{n-1}$ for $n - 1 > k > 0$
$a_0^n = a_1^{n-1}$

**Property 3**

Le $G$ be the function, of two variables $x$ and $y$, defined by the relation $G(x, y) = e^{2xy - y^2}$. We have the relation

$$G(x, y) = \sum_n H_n(x) \frac{y^n}{n!} \quad (A.3)$$

$G(x, y)$ is the generating function of the Hermite polynomials

**Property 4**

For all positive integers $n$ and $m$, one has the relation of orthogonalization

$$\int_{-\infty}^{+\infty} H_n(x) H_m(x) e^{-x^2} dx = 2^n n! \sqrt{\pi} \delta_{nm} \quad (A.4)$$

If we make the variable change

$$x = \frac{t - T}{\sqrt{2}\Delta t} \Rightarrow dx = \frac{dt}{\sqrt{2}\Delta t}$$

we may deduce from the relation (A.4) the relation of orthonormalization

$$\int_{-\infty}^{+\infty} \frac{H_n\left(\frac{t-T}{\sqrt{2}\Delta t}\right) H_m\left(\frac{t-T}{\sqrt{2}\Delta t}\right)}{\sqrt{2^{n+m} n! \, m!}\sqrt{2\pi}\Delta t} e^{-\frac{1}{2}\left(\frac{t-T}{\Delta t}\right)^2} dt = \delta_{nm} \quad (A.5)$$

**Property 5**

$$\int_{-\infty}^{+\infty} t \frac{\left|H_n\left(\frac{t-T}{\sqrt{2}\Delta t}\right)\right|^2}{2^n n! \sqrt{2\pi}\Delta t} e^{-\frac{1}{2}\left(\frac{t-T}{\Delta t}\right)^2} dt = T \quad (A.6)$$

$$\int_{-\infty}^{+\infty} (t-T)^2 \frac{\left|H_n\left(\frac{t-T}{\sqrt{2}\Delta t}\right)\right|^2}{2^n n! \sqrt{2\pi}\Delta t} e^{-\frac{1}{2}\left(\frac{t-T}{\Delta t}\right)^2} dt$$

$$= (2n+1)(\Delta t)^2 = 2\left(n + \frac{1}{2}\right)(\Delta t)^2 \quad (A.7)$$

# Appendix B

**Fourier Transform of a harmonic Gaussian function**

Let us consider a harmonic Gaussian function $\varphi_n$ as defined in the section 2:

$$\varphi_n(t, T, \Omega, \Delta t) = \frac{H_n\left(\frac{t-T}{\sqrt{2}\Delta t}\right)}{\sqrt{2^n n! \sqrt{2\pi}\Delta t}} e^{-\left(\frac{t-T}{2\Delta t}\right)^2} e^{i\Omega t}$$

Our purpose is to prove that the expression of the Fourier transform $\tilde{\varphi}_n$ of the function $\varphi_n$ is

$$\tilde{\varphi}_n(\omega, T, \Omega, \Delta\omega) = \frac{1}{\sqrt{2\pi}} \int_{-\infty}^{+\infty} \varphi_n(t, T, \Omega, \Delta t) e^{-i\omega t} dt$$

$$= \frac{(i)^n H_n\left(\frac{\omega - \Omega}{\sqrt{2}\Delta\omega}\right) e^{-\left(\frac{\omega-\Omega}{2\Delta\omega}\right)^2}}{\sqrt{2^n n! \sqrt{2\pi}\Delta\omega}} e^{-iT(\omega-\Omega)} (B.1)$$

in which $\Delta\omega$ is related to $\Delta t$ by the relation $\Delta t \Delta\omega = \frac{1}{2}$.

Let us perform the calculation

$$\tilde{\varphi}_n(\omega, T, \Omega, \Delta\omega) = \frac{1}{\sqrt{2\pi}} \int_{-\infty}^{+\infty} \varphi_n(t, T, \Omega, \Delta t) e^{-i\omega t} dt$$

$$= \frac{1}{\sqrt{2^n n! \, 2\pi\sqrt{2\pi}\Delta t}} \int_{-\infty}^{+\infty} H_n\left(\frac{t-T}{\sqrt{2}\Delta t}\right) e^{-\left(\frac{t-T}{2\Delta t}\right)^2} e^{-i(\omega-\Omega)t} dt$$

We make the variable change

$$x = \frac{t-T}{\sqrt{2}\Delta t} \Rightarrow dt = \sqrt{2}\Delta t\, dx$$

$\tilde{\varphi}_n(\omega, T, \Omega, \Delta\omega)$

$$= \frac{e^{-iT(\omega-\Omega)}\sqrt{\Delta t}}{\sqrt{2^n n!\,\pi\sqrt{2\pi}}} \int_{-\infty}^{+\infty} H_n(x) e^{-\frac{x^2}{2}-i\sqrt{2}\Delta t(\omega-\Omega)x}\, dx$$

For $a > 0$, let $I_n(a,b)$ be the integral

$$I_n(a,b) = \int_{-\infty}^{+\infty} H_n(x)\, e^{-(ax^2+ibx)} dx \quad (B.2)$$

$$\tilde{\varphi}_n(\omega, T, \Omega, \Delta\omega) = \frac{e^{-iT(\omega-\Omega)}\sqrt{\Delta t}}{\sqrt{2^n n!\,\pi\sqrt{2\pi}}} I_n\left(\frac{1}{2}, \sqrt{2}\Delta t(\omega-\Omega)\right) \quad (B.3)$$

Let us determine the expression of the integral $I_n(a,b)$. We introduce in the calculation the generating function of Hermite polynomials

$$G(x,y) = e^{2xy-y^2} = \sum_n H_n(x)\frac{y^n}{n!}$$

On one side

$$\int_{-\infty}^{+\infty} G(x,y)e^{-(ax^2+ibx)}\, dx = \int_{-\infty}^{+\infty} e^{-ax^2+(2y-ib)x-y^2}\, dx$$

$$= \sqrt{\frac{\pi}{a}} e^{\frac{(2y-ib)^2}{4a}-y^2} = \sqrt{\frac{\pi}{a}} e^{-\frac{b^2}{4a}} e^{-\frac{iby}{a}+\left(\frac{1-a}{a}\right)y^2}$$

$$= \sqrt{\frac{\pi}{a}} e^{-\frac{b^2}{4a}} e^{2\left[\frac{b}{2a}\left(\frac{a}{1-a}\right)^{1/2}\right]\left[i\left(\frac{1-a}{a}\right)^{1/2}y\right]-\left[i\left(\frac{1-a}{a}\right)^{1/2}y\right]^2}$$

$$= \sqrt{\frac{\pi}{a}} e^{-\frac{b^2}{4a}} G\left[\frac{b}{2a}\left(\frac{a}{1-a}\right)^{1/2}, i\left(\frac{1-a}{a}\right)^{1/2}y\right]$$

$$= \sqrt{\frac{\pi}{a}} e^{-\frac{b^2}{4a}} \sum_n H_n\left[\frac{b}{2a}\left(\frac{a}{1-a}\right)^{1/2}\right] i^n \left(\frac{1-a}{a}\right)^{n/2} \frac{y^n}{n!}$$

$$= \sum_n \left[i^n \sqrt{\frac{\pi}{a}} \left(\frac{1-a}{a}\right)^{n/2} H_n\left(\frac{b}{2\sqrt{a(1-a)}}\right) e^{-\frac{b^2}{4a}}\right] \frac{y^n}{n!} \quad (B.4)$$

On the other side

$$\int_{-\infty}^{+\infty} G(x,y)e^{-(ax^2+ibx)}\, dx$$

$$= \sum_n \left[\int_{-\infty}^{+\infty} H_n(x)\, e^{-(ax^2+ibx)} dx\right] \frac{y^n}{n!} \quad (B.5)$$

By comparing the relations (B.4) and (B.5), we may identify

$$I_n(a,b) = \int_{-\infty}^{+\infty} H_n(x)\, e^{-(ax^2+ibx)} dx$$

$$= (i)^n \sqrt{\frac{\pi}{a}} \left(\frac{1-a}{a}\right)^{\frac{n}{2}} H_n\left(\frac{b}{2\sqrt{a(1-a)}}\right) e^{-\frac{b^2}{4a}} \quad (B.6)$$

Then we have for $a = \frac{1}{2}$ and $b = \sqrt{2}\Delta t(\omega-\Omega)$

$$I_n\left(\frac{1}{2}, \sqrt{2}\Delta t(\omega-\Omega)\right)$$

$$= (i)^n \sqrt{2\pi} H_n[\sqrt{2}\Delta t(\omega-\Omega)] e^{-\Delta t^2(\omega-\Omega)^2}$$

Introducing the quantity $\Delta\omega$

$$\Delta\omega\Delta t = \frac{1}{2} \Leftrightarrow \Delta t = \frac{1}{2\Delta\omega} \Leftrightarrow (\Delta t)^2 = \frac{1}{4(\Delta\omega)^2} \quad (B.7)$$

we obtain

$$I_n\left(\frac{1}{2}, \sqrt{2}\Delta t(\omega-\Omega)\right) = (i)^n \sqrt{2\pi} H_n\left(\frac{\omega-\Omega}{\sqrt{2}\Delta\omega}\right) e^{-\left(\frac{\omega-\Omega}{2\Delta\omega}\right)^2}$$

Introducing the relations (B.7) and (B.8) in the relation (B.3) and rearranging, we obtain as expected the relation

$$\tilde{\varphi}_n(\omega, T, \Omega, \Delta\omega) = \frac{(i)^n H_n\left(\frac{\omega-\Omega}{\sqrt{2}\Delta\omega}\right) e^{-\left(\frac{\omega-\Omega}{2\Delta\omega}\right)^2}}{\sqrt{2^n n!\,\sqrt{2\pi}\Delta\omega}} e^{-iT(\omega-\Omega)}$$